\documentclass[runningheads,a4paper]{llncs}

\usepackage{amssymb}
\usepackage[latin1]{inputenc}
\usepackage[T1]{fontenc}

\usepackage{graphicx}

\usepackage{booktabs}

\usepackage{graphicx}
\usepackage{color}
\usepackage[latin1]{inputenc}  
\usepackage{amsmath}           
 
\usepackage{xfrac} 
 
\usepackage{float} 

\usepackage[]{graphicx} 
\usepackage{tikz}

\usetikzlibrary{matrix}

\usepackage{verbatim}
\usetikzlibrary{shapes,arrows,calc}
\usetikzlibrary{arrows}

\usetikzlibrary{arrows,calc}

\begin{document}

\mainmatter  

\title{Cryptosystems Using Automorphisms of Finitely Generated Free Groups}


%
%
\author{Anja I. S. Moldenhauer%
 \and Gerhard Rosenberger}
%

\institute{Fachbereich Mathematik, Universit\"at Hamburg, \\
Bundesstrasse 55, 20146 Hamburg, Germany\\
\texttt{anja.moldenhauer@uni-hamburg.de, gerhard.rosenberger@math.uni-hamburg.de}}

%
%



\maketitle

\begin{abstract}
This paper introduces a newly developed private key cryptosystem and a public key cryptosystem. In the first one, each letter is encrypted with a different key. Therefore, it is a kind of a  one-time pad. The second one is inspired by the ElGamal cryptosystem. Both presented cryptosystems are based on automorphisms of free groups. Given a free group $F$ of finite rank, the automorphism group $Aut(F)$ can be generated by Nielsen transformations, which are the basis of a  linear technique to study free groups and general infinite groups. Therefore Nielsen transformations are introduced. 

\end{abstract}

\textit{Dedicated to Gabriele Kern-Isberner on the occasion of her 60th birthday\footnote{This article appeared 2016 in College Publications Tributes  Volume 29: \textbf{Computational Models of Rationality} \textit{Essays Dedicated to Gabriele Kern-Isberner on the Occasion of Her 60th Birthday}. }.}

\section{Introduction}

The topic of this paper is established in the area of mathematical cryptology, more precisely in group based cryptology.  We refer to the books \cite{BFKR}, \cite{VS} and \cite{MSU} for the interested reader. The books \cite{BFKR} and \cite{VS} can also be used for a first access to the wide area of cryptology.  

We introduce two cryptosystems, the first one is a private key cryptosystem (one-time pad) and the second one is a public key cryptosystem. We require that the reader is familiar with the  general concept of these types of protocols. In cryptology it is common to call the two parties who want to communicate privately with each other Alice and Bob.\\
Throughout the paper let $F$ always be a free group $F=\langle X \mid \phantom{R} \rangle$ of finite rank.
 Both cryptosystems  are based on free groups $F$ of finite rank and automorphisms of $F$. It is known that the group of all automorphisms of $F$, $Aut(F)$, can be generated by Nielsen transformations (see \cite{CRR}).   

 We first review some basic definitions concerning regular Nielsen transformations and Nielsen reduced sets and we give additional information which is also important for the understanding of the paper.\\
Both cryptosystems use automorphisms of a free group $F$ of finite rank. Thus, a random choice of these automorphisms is practical. An approach for this random choice  using the Whitehead-Automorphisms is given. Therefore, the Whitehead-Automorphisms are reviewed. Note, that the Whitehead-Automor\-phisms generate the Nielsen transformations and vice versa, but the use of the Whitehead-Automorphisms is more practical if a random choice of automorphisms is required. \\After this,  a private key cryptosystem using Nielsen transformations and Nielsen reduced sets is introduced. An example  and a security analysis for this private key cryptosystem is given. Finally we explain  a public key cryptosystem which is inspired by the ElGamal cryptosystem, from which we describe two variations and give an example.

The new cryptographic protocols are in part in the dissertation  \cite{M} of A.~\mbox{Moldenhauer}  under her supervisor G. Rosenberger at the University of Hamburg.

\section{Preliminaries for Automorphisms of Free Groups}

 We now review some basic definitions concerning regular Nielsen transformations and Nielsen reduced sets  and we give additional information which will be used later on  (see also \cite{CRR}, \cite{LS}~or~\cite{MKS}).

Let $F$ be a free group on the free generating set $X:=\{x_1,x_2,\ldots, x_q\}$ and let $U:=\{u_1,u_2,\ldots, u_t \} \subset F$, $q,t \geq 2$.  A \textbf{freely reduced word} in $X$ is a word in which the symbols $x_i^{\epsilon}$, $x_i^{-\epsilon}$, for $\epsilon= \pm 1$ and $i=1,2,\ldots,q$, do not occur consecutively. We call $q$ the \textbf{rank of $F$}. The free generating set $X$ is also called a \textbf{basis} of $F$. The elements $u_i$ are freely reduced words with letters in $X^{\pm 1}:=X \cup X^{-1}$, with $X^{-1}:=\{x^{-1}_1,x^{-1}_2,\ldots,x^{-1}_q \}$.

\begin{definition}\label{Nreduceddef}~\\
An \textbf{elementary Nielsen transformation} on \\ $U=\{u_1,u_2,\ldots, u_t \} \subset F$ is one of the following transformations
\begin{itemize}
\item[(T1)] replace some $u_i$ by $u_i^{-1}$;
\item[(T2)] replace some $u_i$ by $u_iu_j$ where $j \ne i$;
\item[(T3)] delete some $u_i$ where $u_i=1$.
\end{itemize}
In all three cases the $u_k$ for $k \ne i$ are not changed.
A (finite) product of elementary Nielsen transformations is called a \textbf{Nielsen transformation}.
A Nielsen transformation is called \textbf{regular} if it is a finite product of the transformations (T1) and (T2), otherwise it is called \textbf{singular}.
The regular Nielsen transformations generate a  group.
The set $U$ is called \textbf{Nielsen-equivalent} to the set $V$, if there is a regular Nielsen transformation from $U$ to $V$. Nielsen-equivalent sets $U$ and $V$ generate the same group, that is, $\langle U \rangle=\langle V \rangle$.
\end{definition}

Now, we agree on some notations. We write $(T1)_i$ if we replace $u_i$ by $u_i^{-1}$ and we write $(T2)_{i.j}$ if we replace $u_i$ by $u_iu_j$. If we want to apply   the same Nielsen transformation $(T2)$ consecutively $t$-times  we write $[(T2)_{i.j}]^t$ and hence replace $u_i$ by $u_iu_j^t$. In all cases the $u_k$ for $k\ne i$ are not changed.

\begin{definition}~\\
A finite set $U$ in $F$ is called \textbf{Nielsen reduced}, if for any three elements $v_1,v_2,v_3$ from $U^{\pm 1}$  the following conditions hold:
\begin{enumerate}
\item[(N0)] $v_1 \ne 1$;
\item[(N1)] $v_1v_2 \ne 1$ implies $|v_1v_2|\geq |v_1|, |v_2|$; 
\item[(N2)] $v_1v_2\ne 1$ and $v_2v_3 \ne 1$ implies $|v_1v_2v_3|>|v_1| - |v_2| + |v_3|$.
\end{enumerate}
Here $|v|$ denotes the \textbf{free length} of $v \in F$, that is, the number of letters from $X^{\pm 1}$ in the freely reduced word $v$.
\end{definition}

\begin{remark}
We say that any word $w$ with finitely many letters from $X^{\pm 1}$ has \textbf{length} $L$ if the number of letters occurring is $L$. The length of a word $w$ is greater than  or equal to the free length of the word $w$. For freely reduced words the length and the free length are equal. If a word $w$  is  not freely reduced then  the length is greater than the free length of $w$. 
\end{remark}

\begin{proposition}\label{propNred}
\textnormal{\cite[Theorem 2.3]{CRR} or \cite[Proposition 2.2]{LS}}\\If $U=\{u_1,u_2,\ldots,u_m\}$ is finite, then $U$ can be carried by a Nielsen transformation into some $V$ such that $V$ is Nielsen reduced. We have $\text{rank} (\langle V \rangle) \leq m$.
\end{proposition}

\begin{proposition}\textnormal{\cite[Corollary 3.1]{MKS}}\\
Let $H$ be a finitely generated subgroup of the free group $F$ on the free generating set $X$. Let $U=\{u_1,u_2,\ldots,u_t\}$, $u_i$ words in  $X$, be a Nielsen reduced set. Then, out of all systems of generators for $H$, the set $U$ has the shortest total $x$-length, that is $\sum_{i=1}^{t} |u_i|$.
\end{proposition}

\begin{remark}
If $F_V$ is a finitely generated subgroup of $F=\langle X \mid \phantom{R} \rangle$, with free generating set $V=\{v_1,v_2,\ldots,v_N\}$, $v_i$ words in  $X$, then there exist only finitely many Nielsen reduced sets $U_i=\{u_{i_1}, u_{i_2},\ldots, u_{i_N}\}$, $i=1,2,\ldots,\ell$, to $V$, which are Nielsen-equivalent. With the help of a lexicographical order $<_{lex}$ (see for instance \cite[Proof of Satz 2.3]{CRR}) the smallest set $U_{s}$, in the set of all Nielsen reduced sets $U^V_{Nred}:=\{U_1,U_2,\ldots,U_{\ell}\}$ to $V$, can be uniquely marked. With the use of regular Nielsen transformations it is possible to obtain this marked  set $U_{s}$ starting from any arbitrary set in $U^V_{Nred}$.
\end{remark}

\begin{proposition}\label{NTgenAutF} \textnormal{\cite[Korollar 2.10]{CRR}}\\ 
Let $F$ be the free group of rank $q$. Then, the group of all automorphisms of $F$, Aut(F), is generated by the elementary Nielsen transformations (T1) and (T2).\\
More precisely: Each automorphism of $F$ is describable as a regular Nielsen transformation between two bases of $F$, and, each regular Nielsen transformation between two bases of $F$ defines an automorphism of $F$.

\end{proposition}

\begin{remark}\label{AlgNred}
In \cite{S} an  algorithm, using elementary Nielsen transformations, is presented which, given a finite set $S$ of $m$ words of a free group, returns a set $S'$ of Nielsen reduced words such that $\langle S \rangle = \langle S' \rangle$; the algorithm runs in $\mathcal O (\ell^2 m^2)$ time, where $\ell$ is the maximum length of a word in $S$.
\end{remark}

\begin{theorem}\label{NredFrei}
\textnormal{\cite[Satz 2.6]{CRR}}\\
Let $U$ be Nielsen reduced, then $\langle  U \rangle$ is free on $U$.
\end{theorem}

For the next lemma we need some notations.
Let $w\ne 1$ be a freely reduced word in $X$. The initial segment $s$ of $w$ which is ``a little more than half'' of $w$ (that is, $\frac{1}{2}|w|< |s| \leq \frac{1}{2}|w|+1$) is called the \textbf{major initial segment} of $w$. The \textbf{minor initial segment} of $w$ is that initial segment $s'$ which is ``a little less than half'' of $w$ (that is, $\frac{1}{2}|w|-1 \leq |s'| < \frac{1}{2} |w|$). Similarly, \textbf{major} and \textbf{minor terminal segments} are defined.\\
If the free length of the word $w$ is even, we call the initial segment $s$ of $w$, with  $|s|=\frac{1}{2}|w|$ the \textbf{left half} of $w$. Analogously, we call the terminal segment $s'$ of $w$ with $|s'|=\frac{1}{2}|w|$ the \textbf{right half} of $w$.\\
Let $\{w_1,w_2, \ldots, w_n\}$ be a set of freely reduced words in $X$, which are not the identity. An initial segment of a $w$-symbol (that is, of either $w_i$ or $w^{-1}_i$, which are different $w$-symbols) is called \textbf{isolated} if it does not occur as an initial segment of any other $w$-symbol. Similarly, a terminal segment is isolated if it is a terminal segment of a unique $w$-symbol.

\begin{lemma}\label{NTreducleicht}
\textnormal{\cite[Lemma 3.1]{MKS}}\\
Let $M:=\{w_1,w_2,\ldots,w_m\}$ be a set of freely reduced words in $X$ with $w_j\ne 1$, $1 \leq j \leq m$. Then $M$ is Nielsen reduced if and only if the following conditions are satisfied:
\begin{enumerate}
\item Both the major initial and major terminal segments of each $w_i \in M$ are isolated.
\item For each $w_i\in M$ of even free length, either its left half or its right half is isolated.
\end{enumerate}
\end{lemma}

\begin{definition}\label{defprim}~\\
Let $F$ be a free group of rank $q$ and let $G$ be a free subgroup of $F$ with rank $m$. An element $g \in G$ is called a \textbf{primitive element} of $G$, if  a basis $U$ of $G$ with $g \in U$ exists.
\end{definition}

\begin{proposition}\label{primnumb1}
 \textnormal{\cite{MS}}\\
The number of primitive elements of free length $k$ of the free group\\ $F_2 = \langle x_1,x_2 \mid \phantom{R} \rangle$ (and therefore, in any group $F_q=\langle x_1,x_2,\ldots, x_q \mid \phantom{R}\rangle$, $q \geq 2$) is:
\begin{enumerate}
\item more than $\frac{8}{3\sqrt{3}}\cdot (\sqrt{3})^k$ if $k$ is odd;
\item more than $\frac{4}{3} \cdot (\sqrt{3})^k$ if $k$ is even.
\end{enumerate}
\end{proposition}

\begin{theorem}\label{primnumb2}
\textnormal{\cite{BMS}}\\
If $P(q,k)$ is the number of primitive elements of free length $k$ of the free group $F_q=\langle x_1,x_2,\ldots, x_q \mid \phantom{R}\rangle$, $q \geq 3$, then for some constants $c_1$, $c_2$, we  have
$$
c_1 \cdot (2q-3)^k \leq P(q,k) \leq c_2 \cdot (2q-2)^k.
$$
\end{theorem}

\begin{definition}\textnormal{\cite{R}}\\
A subgroup $H$ of $F$ is called \textbf{characteristic} in $F$ if $\varphi(H)=H$ for every automorphism $\varphi$ of $F$. 
\end{definition}

For $ n \in \mathbb N$ let $\mathbb Z_n:= \sfrac{\mathbb Z}{ n \mathbb  Z}$ be the ring of integers modulo $n$. The corresponding residue class in $\mathbb Z_n$ for an integer $\beta$ is denoted by $\overline{\beta}$ (see also \cite{BFKR}). 

\begin{definition}\textnormal{\cite{BFKR}}\label{defLCG}\\
Let $n \in \mathbb N$ and $ \overline \beta, \overline \gamma \in \mathbb Z_n$. A bijective mapping $h: \mathbb Z_n \to \mathbb Z_n$ given by $\mathrm x \mapsto \overline{\beta} \mathrm x + \overline \gamma$ is called a \textbf{linear congruence generator}.
\end{definition}

\begin{theorem}\textnormal{\cite{BFKR} (Maximal period length for $n=2^m$, $m\in \mathbb N$)}\label{LCGML}\\
Let $n \in \mathbb N$, with $n=2^{m}$, $m \geq 1$ and let $\beta, \gamma \in \mathbb Z$ such that $h: \mathbb Z_n~\to~\mathbb Z_n$, with $\mathrm x \mapsto \overline \beta \mathrm x + \overline \gamma$, is a linear congruence generator. Further let $\alpha \in \{0,1,\ldots,{n-1}\}$ be given and $\mathrm x_1= \overline \alpha$, $\mathrm x_2=h(\mathrm x_1)$, $\mathrm x_3=h(\mathrm x_2)$, $\ldots$.\\ Then the sequence $\mathrm x_1,\mathrm x_2,\mathrm x_3,\ldots$ is periodic with maximal periodic length $n=2^m$ if and only if the following holds:
\begin{enumerate}
\item $\beta$ is odd.
\item If $m \geq 2$ then $\beta \equiv 1 \pmod 4$.
\item $\gamma$ is odd.
\end{enumerate}
\end{theorem}

\begin{theorem}\label{TL}\textnormal{\cite{L}}\\
Let $F$ be a free group with countable number of generators $x_1,x_2,\ldots$. Corresponding to $x_j$ define 
\begin{align*}
M_j=\begin{pmatrix}
 -r_j & \ \ -1+r_j^2\\
   1  & \ \ -r_j
\end{pmatrix}
\end{align*}
with $r_j \in \mathbb Q$ and 
\begin{align*}
r_{j+1}-r_j &\geq 3\\ 
r_1 &\geq 2.
 \end{align*}
Then G* generated by $\{M_1,M_2,\ldots\}$ is isomorphic to $F$. 
\end{theorem}

\section{The Random Choice of the Automorphisms of $Aut(F)$}\label{chooseNT}

 Let $F=\langle X \mid  \phantom{R} \rangle$ be the free group on the free generating set $X$ with  $|X|~=~q$. The cryptosystems we develop are based on automorphisms of  $F$. These automorphisms should be chosen randomly. It is known, see  Proposition \ref{NTgenAutF}, that the Nielsen transformations generate the automorphism group $Aut(F)$. For a realization of a random choice procedure the Whitehead-Automorphisms will be used. 
 
 A fixed set of randomly chosen automorphisms is part of the key space for the private key cryptosystem.

\begin{definition}
Whitehead-Automorphisms:

\begin{enumerate}
\item Invert the letter $a$ and leave all other letters invariant: 
\begin{align*}
i_a(b) = \begin{cases} a^{-1} &\text{ for } a=b \\ 
                       b      & \text{ for } b \in X \setminus\{a\}.\end{cases}
\end{align*}
There are $q$ \textbf{Whitehead-Automorphisms} of this type.

\item Let $a \in X$ and $L,R,M$ be three pairwise disjoint subsets of $X$, with $a \in M$. Then the tuple  $(a,L,R,M)$ defines a  \textbf{Whitehead-Automorphisms}\\ $W_{(a,L,R,M)}$ as follows
\begin{align*}
W_{(a,L,R,M)}(b) = \begin{cases} ab &\text{ for } b \in L \\ 
ba^{-1} &\text{ for } b \in R \\
								 aba^{-1} &\text{ for } b \in M \\
                       b      & \text{ for } b \in X \setminus ( L \cup M \cup R ). \end{cases}
\end{align*}
There are $q \cdot 4^{q-1}$ automorphisms of this type. 
\end{enumerate}
\end{definition}

Note, that $W^{-1}_{(a,L,R,M)}= i_a \circ W_{(a,L,R,M)} \circ i_a$.\\

With this definition it is clear how the Whitehead-Automorphisms can be generated as a product of regular Nielsen transformations. Conversely, the Whitehead-Automorphisms generate the group of the Nielsen transformations and therefore also the automorphism group $Aut(F)$ (see also \cite{DKR}). With the Whitehead-Automorphisms it is simple to realize a random choice of automorphisms. We now give an approach for this choice.\\

\underline{\textbf{An approach for choosing randomly automorphisms of $Aut(F)$:}}\\

Let $X=\{x_1,x_2,\ldots,x_q\}$ be the free generating set for the free group $F$.
\begin{enumerate}
\item First of all it should be decided in which order an automorphism $f_i$ is generated by  automorphisms of type $i_a$ and $W_{(a,L,R,M)}$.  For this purpose an automorphism of type $i_a$ is identified with  a zero and $W_{(a,L,R,M)}$ with a  one. A sequence of zeros and ones is randomly generated. This sequence is translated to randomly chosen Whitehead-Automorphisms and hence presents an automorphism $f_i \in Aut(F)$. This translation is as follows:

 \item[2.1.] For a zero in the sequence we generate $i_a$ randomly: choose a random number $z$, with $1\leq z \leq q$; hence an element $a \in X$ must be chosen to declare the automorphism. Then it is $a:=x_{z}$ and   hence $x_z$ is replaced by $x^{-1}_{z}$ and all other letters are invariant.

 \item[2.2.] For a one in the sequence we generate $W_{(a,L,R,M)}$ randomly: choose a random number  $z$, with $1\leq z \leq q$. Hence it is $a:=x_z$. Moreover it is $a \in M$. After this  the disjoint sets $L,R,M$ $\subset X$ are chosen randomly. One possible approach is the following:

\begin{enumerate}
\item Choose random numbers $z_1, z_2$ and $z_3$ with 
\begin{align*}
0 &\leq z_1 \leq q-1, \text{}\\
0 &\leq z_2 \leq q-1-z_1,\\
0 &\leq z_3 \leq q-1-z_1-z_2.
\end{align*}

If we are in the situation of $z_1=z_2=z_3=0$ we get the identity $id_{X}$. If this case arises a random number $\tilde{z}$ from the set $\{1,2,\ldots,q\} \setminus \{z\}$ is chosen and hence the element $x_{\tilde{z}}$ is assigned randomly to one of the sets $L$, $R$ or $M$;  therefore the identity is avoided.

It is 
\begin{align*}
|L| = z_1, \qquad |R|= z_2, \qquad |M|=z_3 + 1.
\end{align*}
\item Choose $z_1$  pairwise different random numbers $\{r_1,r_2,\ldots,r_{z_1}\}$ of the set 
$\{1,2,\ldots,q\} \setminus \{z\}$. Then $L$ is the set  
\begin{align*}
L=\{ x_{r_1}, x_{r_2},\ldots, x_{r_{z_1}}\}.
\end{align*}
\item Choose $z_2$  pairwise different random numbers $\{p_1,p_2,\ldots,p_{z_2}\}$ of the set 
$\{1,2,\ldots,q\} \setminus \left( \{z\} \cup \{r_1,r_2,\ldots,r_{z_1}\} \right)$. Then  $R$ is the set
\begin{align*}
R=\{ x_{p_1}, x_{p_2},\ldots, x_{p_{z_2}}\}.
\end{align*}
\item Choose $z_3$  pairwise different random numbers $\{t_1,t_2,\ldots,t_{z_3}\}$  of the set 
$\{1,2,\ldots,q\} \setminus  \left( \{z\} \cup \{r_1,r_2,\ldots,r_{z_1}\} \cup \{p_1,p_2,\ldots,p_{z_2}\} \right)$. Then  $M$ is the set 
\begin{align*}
M=\{ x_{t_1}, x_{t_2},\ldots, x_{t_{z_3}}\} \cup \{ a \}.
\end{align*}
\end{enumerate}
\end{enumerate}

\begin{remark}
If Alice and Bob use Whitehead-Automorphisms to generate automorphisms on a free group with free generating set $X$ they should take care, that there are no sequences of the form
\begin{enumerate}
\item $i_a \circ i_a=id_{X}$,
\item $W_{(a,L,R,M)} \circ \underbrace{i_a \circ W_{(a,L,R,M)} \circ i_a}_{=W^{-1}_{(a,L,R,M)}}= id_{X}$ \quad or \\  $  \underbrace{i_a \circ W_{(a,L,R,M)} \circ i_a}_{=W^{-1}_{(a,L,R,M)}} \circ
   W_{(a,L,R,M)}=id_{X},$
\end{enumerate}
for the automorphism $f_j$.
They also should not use Whitehead-Auto\-morphisms sequences for $f_j$, which cancel each other and so be vacuous for the encryption.
\end{remark}

\section{Private Key Cryptosystem Based on Automorphisms of Free Groups $F$}\label{privkeycryp1}

Before Alice and Bob are able to communicate with each other, they have to make some arrangements. 

\begin{center}\underline{Public Parameters}\end{center}
They first agree on the public parameters.
\begin{enumerate}
\item A free group $F$ with free generating set $X=\{x_1,x_2,\ldots,x_q\}$, with $q \geq 2$.
\item A plaintext alphabet $A=\{a_1,a_2,\ldots,a_N\}$, with $N \geq 2$.
\item A subset $\mathcal F_{aut}:=\{f_1,f_2,\ldots,f_{2^{128}}\} \subset Aut(F)$  of automorphisms of $F$ is chosen. It is $f_i~:~F \to F$ and the $f_i$, $i=1,2,\ldots,2^{128}$,  pairwise different, are generated with the help of $0$-$1$-sequences (of different length) and random numbers as described in Section~\ref{chooseNT}. The set   $\mathcal F_{aut}$ is part of the key space.
\item They agree on a linear congruence generator $h:\mathbb Z_{2^{128}} \to \mathbb Z_{2^{128}}$ with a maximal period length (see Definition \ref{defLCG} and Theorem \ref{LCGML}). 
\end{enumerate}

\begin{remark}
If the set $\mathcal F_{aut}$ and the linear congruence generator $h$ are public Alice and Bob are able to change the automorphisms and the generator publicly  without a private meeting. 
The set  $\mathcal F_{aut}$ should be large enough to make a brute force search ineffective.\\
Another variation could be, that Alice and Bob choose the number of elements in the starting set $\mathcal F_{aut}$  smaller than $2^{128}$, say for example  $2^{10}$. These starting automorphism set $\mathcal F_{aut}$ should be chosen privately by Alice and Bob as their set of seeds  and should not be made public. Then Alice and Bob can extend publicly the starting set $\mathcal F_{aut}$ to the set $\mathcal F_{aut_1}$ of automorphisms such that $\mathcal F_{aut_1}$ contains, say for example, $2^{32}$ automorphisms. The number of all elements in $\mathcal F_{aut_1}$  should make a brute force attack inefficient. The linear congruence generator stays analogously, just the domain and codomain must be adapted to, say for example, $\mathbb Z_{2^{32}}$. Because of  Theorem \ref{LCGML} Alice and Bob get at all times a linear congruence generator with maximal periodic length.
\end{remark}


\begin{center}\underline{Private  Parameters}\end{center}
Now they agree on the private parameters.

\begin{enumerate}
\item A free subgroup $F_U$ of $F$ with rank $N$ and the free generating set\\ $U=\{u_1,u_2,\ldots,u_{N}\}$ is chosen where $U$ is a minimal  Nielsen reduced set (with respect to a lexicographical order) and the $u_i$ freely reduced words in $X$.
 Such systems $U$ are easily to construct using Theorem~\ref{NredFrei} and Lemma~\ref{NTreducleicht} (see also \cite{CRR} and \cite{LS}). It is $\mathcal U_{Nred}$ the set of all minimal Nielsen reduced sets with $N$ elements in $F$, which is part of the key space.
 \item They use  a one to one  correspondence
\begin{align*}
A &\to U\\
a_j &\mapsto u_j \qquad \text{ for } j=1,\ldots,N.
\end{align*}

\item Alice and Bob agree on an automorphism $f_{\overline{\alpha}} \in \mathcal F_{Aut}$, hence  $\alpha$ is the common secret starting point $\alpha \in \{0,1,\ldots,2^{128}-1\}$, with $\mathrm x_1= \overline \alpha \in \mathbb Z_{2^{128}}$, for the linear congruence generator. With this $\alpha$ they are able to  generate the sequence of automorphisms of the set $\mathcal F_{aut}$, which they use for encryption and decryption, respectively.
\end{enumerate}

\textbf{The key space:} The set $\mathcal U_{Nred}$ of all minimal (with respect to a lexicographical order) Nielsen reduced subsets of $F$ with $N$ elements. The set $\mathcal F_{Aut}$ of $2^{128}$ randomly chosen automorphisms of $F$.

\begin{center}\underline{Protocol}\end{center}

Now we explain the protocol and look carefully at  the steps for Alice and Bob.\\

\textbf{Public knowledge:} $F=\langle X \mid \phantom{R}\rangle$, $X=\{x_1,x_2,\ldots,x_q\}$ with $q \geq 2$; plaintext alphabet $A=\{a_1,a_2,\ldots,a_N\}$ with $N \geq 2$; the set $\mathcal F_{Aut}$; a linear congruence generator $h$.\\

\textbf{Encryption and Decryption Procedure:}  
\begin{enumerate}
\item Alice and Bob agree privately on a set $U \in \mathcal U_{Nred}$ and an automorphism $f_{\overline{\alpha}} \in \mathcal F_{Aut}$.   They also know the one to one correspondence between $U$ and $A$.
\item Alice wants to transmit the message  
$$S=s_1s_2\cdots s_z, \quad  z \geq 1,$$
with $s_i \in A$ to Bob.

\item[2.1.] Alice generates with the linear congruence generator $h$ and the knowledge of $f_{\overline{\alpha}}$ the $z$ automorphisms $f_{\mathrm x_1}, f_{\mathrm x_2}, \ldots, f_{\mathrm x_z}$, which she needs for encryption. It is $\mathrm x_1=\overline{\alpha}, \mathrm x_2=h(\mathrm x_1), \ldots, \mathrm x_{z}=h(\mathrm x_{z-1})$.

\item[2.2.] 
 The encryption is as follows
\begin{align*}
\text{if } s_i=a_t \quad \text{then } s_i \mapsto c_i:= f_{\mathrm x_i}(u_t), \quad 1\leq i \leq z, \quad 1\leq t \leq N.
\end{align*}
Recall that the one to one correspondence $A\to U$ with $a_j \mapsto u_j$, for $j=1,2,\ldots,N$, holds. 
The ciphertext 
\begin{align*}
C&= f_{\mathrm x_1}(s_1)f_{\mathrm x_2}(s_2)\cdots f_{\mathrm x_z}(s_z)\\
 &= c_1c_2\cdots c_z
\end{align*}
is sent to Bob. We call $c_j$ the ciphertext units. We do no cancellations between $c_i$ and $c_{i+1}$, for $1\leq i \leq z-1$.

\item Bob gets the ciphertext 
$$
C=c_1c_2\cdots c_z,
$$ 
and the information that he has to  use $z$ automorphisms of $F$ from the set $\mathcal F_{Aut}$ for decryption. 
He has now two possibilities for decryption.
\item[3.1.a.] With the knowledge of $f_{\overline{\alpha}}$, the linear congruence generator $h$ and the number $z$, he  computes for each automorphism $f_{\mathrm x_i}$, $i=1,2,\ldots,z$,  the inverse automorphism $f^{-1}_{\mathrm x_i}$.
\item[3.1.b.] With the knowledge of $f_{\overline{\alpha}}$, the set $U=\{u_1,u_2,\dots,u_N\}$, the linear congruence generator $h$ and the number $z$,  he  computes for each automorphism~$f_{\mathrm x_i}$, $i=1,2,\ldots,z$,  the set 
$$
U_{f_{\mathrm x_i}}=\{f_{\mathrm x_i}(u_1), f_{\mathrm x_i}(u_2),\ldots, f_{\mathrm x_i}(u_N)\}.
$$
Hence, with the one to one correspondence between $U$ and $A$,  he gets a one to one correspondence between the letters in the alphabet $A$ and the words of the ciphertext depending on the automorphisms $f_{\mathrm x_i}$. This is shown in Table~\ref{ZuordnungCRyptoFabstract}. 

\begin{table}[H]
\caption{Plaintext alphabet $A=\{a_1,a_2,\ldots,a_N\}$ corresponded to ciphertext alphabet $U_{f_{\mathrm x_i}}$ depending on the automorphisms $f_{\mathrm x_i}$ }\label{ZuordnungCRyptoFabstract}
\[
\begin{array}{c||c | c | c | c}

 &   &  &&\\
& \ U_{f_{\mathrm x_1}}  \ & \ U_{f_{\mathrm x_2}}  \ & \quad \cdots \quad & \ U_{f_{\mathrm x_z}}  \\
 &   &  &&\\
\hline
\hline
 &   &  &&\\
a_1 & f_{\mathrm x_1}(u_1)   & f_{\mathrm x_2}(u_1) & \cdots & f_{\mathrm x_z}(u_1)\\
 &   &  &&\\
a_2 & f_{\mathrm x_1}(u_2)  &  f_{\mathrm x_2}(u_2) & \cdots & f_{\mathrm x_z}(u_2)\\
 &   &  &&\\
\vdots & \vdots   & \vdots  & \cdots & \vdots \\
 &   &  &&\\
a_N & f_{\mathrm x_1}(u_N)  &  f_{\mathrm x_2}(u_N) & \cdots & f_{\mathrm x_z}(u_N) \\
 &   &  &&\\
\end{array}
\]
\end{table}

\item[3.2.] With the knowledge of the Table \ref{ZuordnungCRyptoFabstract} or the inverse automorphisms $f^{-1}_{\mathrm x_i}$, respectively, the decryption is as follows
\begin{align*}
\text{if } c_i=f_{\mathrm x_i}(u_t) \quad \text{then } c_i \mapsto s_i:= f^{-1}_{\mathrm x_i}(c_i)=a_t, \quad 1\leq i \leq z, \quad 1\leq t \leq N.
\end{align*}
He generates the plaintext message
\begin{align*}
S&=f^{-1}_{\mathrm x_1}(c_1)f^{-1}_{\mathrm x_2}(c_2)\cdots f^{-1}_{\mathrm x_z}(c_z)\\
&=s_1s_2\cdots s_z,
\end{align*}
with $s_i \in A$, from Alice.
\end{enumerate}

\begin{remark}
The cryptosystem is a polyalphabetic system. A word $u_i \in U$, and hence a letter $a_i \in A$, is encrypted differently at different places in the plaintext.
\end{remark}

\begin{example}\label{example1}
This example was executed with the help of the computer program GAP and the package ``FGA\footnote{Free Group Algorithms, a GAP4 Package by Christian Sievers, TU Braunschweig.}''.

First Alice and Bob agree on the \textbf{public parameters}.
\begin{enumerate}
\item Let $F$ be the free group on the free generating set $X=\{a,b,c,d\}$.
\item 
Let $\tilde A:=\{a_1,a_2,\ldots,a_{12}\}=\{A,E,I,O,U,T,M,L,K,Y,B,N\}$ be the plaintext alphabet.
\item A set $\mathcal F_{Aut}$ is determined. In this example we give the automorphisms, which Alice and Bob use for encryption and decryption, respectively, just at the moment when they are needed.
\item The linear congruence generator with maximal periodic length is 
\begin{align*}
h: \mathbb Z_{2^{128}} &\to \mathbb Z_{2^{128}}\\
                    \mathrm x    &\mapsto  \overline 5 \mathrm x + \overline 3.
\end{align*}
\end{enumerate}

The \textbf{private parameters} for this example are:
\begin{enumerate}
\item The free group $F_{\tilde U}$ of $F$ with the free generating set 
\begin{align*}
\tilde{U}=&\{u_1,u_2,\ldots,u_{12}\}\\
		=&\{ba^2, cd, d^2c^{-2}, a^{-1}b, a^4b^{-1}, b^3a^{-2}, bc^3, bc^{-1}bab^{-1}, \\ & \phantom{P}c^2ba, c^2dab^{-1}, a^{-1}d^3c^{-1}, a^2 d b^2 d^{-1} \}.
\end{align*}
It is known, that $a_i \mapsto u_i$, $i=1,2,\ldots,12$, for $u_i \in \tilde U$ and $a_i \in \tilde A$.
The set $\tilde U$ is a Nielsen reduced set and the group $F_{\tilde U}$ has rank $12$. 
 Alice and Bob agree on the starting automorphism $f_{\overline{93}}$, hence it is  $\mathrm x_1=\overline{\alpha}=\overline{93}$.
\end{enumerate}

We look at the \textbf{encryption and decryption procedure} for Alice and Bob.
\begin{enumerate}

\item[2.] With the above agreements \textbf{Alice} is able to encrypt her message 
\begin{center}
 S = \texttt{I LIKE BOB}.
 \end{center}
Her message is of length $8$. She generates the ciphertext as follows:

\item[2.1.] She first determines, with the help of the linear congruence generator $h$, the automorphisms $f_{\mathrm x_i}$, $i=1, 2, \ldots, 8$, which she needs for encryption. It is
\begin{alignat*}{8}
&\mathrm x_1= \overline \alpha &&=& \overline{93}, \qquad &\mathrm x_2=h(\mathrm x_1)&&=& \overline{468}, \qquad &\mathrm x_3=h(\mathrm x_2)&&=& \overline{2343}, \\
 &\mathrm x_4=h(\mathrm x_3)&&=&\overline{11718},  \qquad &\mathrm x_5=h(\mathrm x_4)&&= &\overline{58593}, \qquad &\mathrm x_6=h(\mathrm x_5)&&=& \overline{292968}, \\ 
&\mathrm x_7=h(\mathrm x_6)&&=&\overline{1464843}, \qquad &\mathrm x_8=h(\mathrm x_7)&&=&\overline{7324218}. \qquad &\phantom{\mathrm x_6=h(\mathrm x_5)}&&\phantom{=}&\phantom{292968} 
\end{alignat*}
The automorphisms are described with the help of regular Nielsen transformations, it is 

\begin{enumerate}
\item[$f_{\mathrm x_1}$=]  $(N1)_3  (N2)_{1.4}   (N2)_{4.3}  (N2)_{2.3}   (N1)_3    (N2)_{1.4}   (N2)_{3.1}$,
\begin{align*}
f_{\mathrm x_1}: F &\to F\\
	a &\mapsto ad^2c^{-1},  \
	b \mapsto  bc^{-1}, \
	c \mapsto  cad^2c^{-1},  \
	d \mapsto  dc^{-1};
\end{align*}
\item[$f_{\mathrm x_2}$=] $(N2)_{1.4}  (N1)_2  (N2)_{2.4}  (N2)_{3.1}  (N1)_2  (N1)_1  (N2)_{1.3}  [(N2)_{4.3}]^2  (N1)_3  $, 
\begin{align*}
f_{\mathrm x_2}: F &\to F\\
	a &\mapsto d^{-1}a^{-1}cad,  \
	b \mapsto  d^{-1}b, \
	c \mapsto  d^{-1}a^{-1}c^{-1}, \
	d \mapsto  d(cad)^2;
\end{align*}
\item[$f_{\mathrm x_3}$=] $(N1)_2   (N2)_{4.2}   (N1)_4   (N2)_{2.4}   (N1)_2   (N2)_{4.2}   (N1)_3  (N2)_{2.1}   (N2)_{3.2}$ \\ $ [(N2)_{1.4}]^3   (N1)_2   (N2)_{4.2}$,
\begin{align*}
f_{\mathrm x_3}: F &\to F\\
	a &\mapsto ab^3,  \
	b \mapsto  a^{-1}d^{-1},\
	c \mapsto  c^{-1}da, \
	d \mapsto  ba^{-1}d^{-1};
\end{align*}
\item[$f_{\mathrm x_4}$=] $[(N2)_{3.1}]^2   (N1)_2  [(N2)_{2.1}]^3  (N2)_{2.4}   (N2)_{4.2}  (N2)_{1.3}$, 
\begin{align*}
f_{\mathrm x_4}: F &\to F\\
	a &\mapsto aca^2,  \
	b \mapsto  b^{-1}a^3d, \
	c \mapsto  ca^2, \
	d \mapsto  db^{-1}a^3d;
\end{align*}
\item[$f_{\mathrm x_5}$=] $(N2)_{1.2} (N1)_3  (N1)_1   [(N2)_{4.3}]^2  (N2)_{1.2}  (N1)_2   (N1)_3  (N2)_{2.4}  (N2)_{3.1}$, 
\begin{align*}
f_{\mathrm x_5}: F &\to F\\
	a &\mapsto b^{-1}a^{-1}b,  \
	b \mapsto  b^{-1}dc^{-2}, \
	c \mapsto  cb^{-1}a^{-1}b, \
	d \mapsto  dc^{-2};
\end{align*}
\item[$f_{\mathrm x_6}$=] $(N1)_1   (N2)_{2.3}    (N2)_{3.1}    (N1)_2   (N2)_{1.2}   (N2)_{4.2}$, 
\begin{align*}
f_{\mathrm x_6}: F &\to F\\
	a &\mapsto a^{-1}c^{-1}b^{-1},   \
	b \mapsto  c^{-1}b^{-1},\
	c \mapsto  ca^{-1}, \
	d \mapsto  dc^{-1}b^{-1};
\end{align*}
\item[$f_{\mathrm x_7}$=] $[(N2)_{2.1}]^3    (N1)_3  [(N2)_{4.3}]^3  (N1)_1  (N2)_{1.2}  (N1)_2   (N2)_{2.4}  (N2)_{3.1}$, 
\begin{align*}
f_{\mathrm x_7}: F &\to F\\
	a &\mapsto a^{-1}ba^{3},  \
	b \mapsto  a^{-3}b^{-1}dc^{-3}, \
	c \mapsto  c^{-1}a^{-1}ba^3, \
	d \mapsto  dc^{-3};
\end{align*}
\item[$f_{\mathrm x_8}$=] $(N2)_{1.4}   (N1)_2  (N1)_3   (N2)_{2.1}   [(N2)_{3.4}]^2  (N1)_4    (N1)_1  (N1)_3  (N2)_{4.2}$, 
\begin{align*}
f_{\mathrm x_8}: F &\to F\\
	a &\mapsto d^{-1}a^{-1},  \
	b \mapsto  b^{-1}ad, \
	c \mapsto  d^{-2}c, \
	d \mapsto  d^{-1}b^{-1}ad.
\end{align*}
Note, that the regular Nielsen transformations are applied from the left to the right.
\end{enumerate}

\item[2.2] The ciphertext is now 
\begin{align*}
C=&f_{\mathrm x_1}(\texttt{I})f_{\mathrm x_2}(\texttt{L})f_{\mathrm x_3}(\texttt{I})f_{\mathrm x_4}(\texttt{K})f_{\mathrm x_5}(\texttt{E})f_{\mathrm x_6}(\texttt{B})f_{\mathrm x_7}(\texttt{O})f_{\mathrm x_8}(\texttt{B})\\
=& f_{\mathrm x_1}(d^2c^{-2})f_{\mathrm x_2}(bc^{-1}bab^{-1})f_{\mathrm x_3}(d^2c^{-2})f_{\mathrm x_4}(c^2ba)f_{\mathrm x_5}(cd)f_{\mathrm x_6}(a^{-1}d^3c^{-1})\\ &f_{\mathrm x_7}(a^{-1}b)f_{\mathrm x_8}(a^{-1}d^3c^{-1})\\
=&dc^{-1}d^{-1}a^{-1}d^{-2}a^{-1}c^{-1}  \wr \  d^{-1}bcabd^{-1}a^{-1}cadb^{-1}d \ \wr \\ 
&   (ba^{-1}d^{-1})^2(a^{-1}d^{-1}c)^2  \wr 
(ca^2)^2b^{-1}a^3daca^2  \wr cb^{-1}a^{-1}bdc^{-2}  \wr \\
&  bca(dc^{-1}b^{-1})^3ac^{-1}  \ \wr  a^{-1} (a^{-2} b^{-1})^2 d c^{-3}  \wr \ (ab^{-1})^3adc^{-1}d^2 \\
=& c_1 c_2 c_3 c_4 c_5 c_6 c_7 c_8.
\end{align*}
The symbol `` $\wr$'' marks the end of a ciphertext unit $c_i$.

\item[3.] \textbf{Bob} gets the ciphertext 
\begin{align*}
C=&dc^{-1}d^{-1}a^{-1}d^{-2}a^{-1}c^{-1}  \wr \  d^{-1}bcabd^{-1}a^{-1}cadb^{-1}d \ \wr \\ 
&   (ba^{-1}d^{-1})^2(a^{-1}d^{-1}c)^2  \wr 
(ca^2)^2b^{-1}a^3daca^2  \wr cb^{-1}a^{-1}bdc^{-2}  \wr \\
&  bca(dc^{-1}b^{-1})^3ac^{-1}  \ \wr  a^{-1} (a^{-2} b^{-1})^2 d c^{-3}  \wr \ (ab^{-1})^3adc^{-1}d^2 
\end{align*}
from Alice. Now he knows, that he needs eight automorphisms for decryption. 
\begin{enumerate}
\item[3.1.] Bob knows the set $U$, the linear congruence generator $h$ and the starting seed automorphism $f_{\overline{93}}$. For decryption he uses tables (analogously to Table~\ref{ZuordnungCRyptoFabstract}).

Now, he is able to compute for each automorphism $f_{\mathrm x_i}$ the set $U_{f_{\mathrm x_i}}$, $i=1,2,\ldots,8$, and  to generate the tables  Table~\ref{ZuordnungCRyptoFI}, Table~\ref{ZuordnungCRyptoFII},  Table~\ref{ZuordnungCRyptoFIII} and Table \ref{ZuordnungCRyptoFIV}.

\begin{table}[H]
\caption{Correspondence: plaintext alphabet to ciphertext alphabet I}\label{ZuordnungCRyptoFI}
\begin{scriptsize}\[
\begin{array}{c||c|c}
 &   &  \\
& \ U_{f_{\mathrm x_1}}  \ & \ U_{f_{\mathrm x_2}}   \ \\
 &   &  \\
\hline
\hline
 &   &  \\
A &  b(c^{-1}ad^2)^2c^{-1} &  d^{-1} b d^{-1} a^{-1} c^2 a d \\
 &   &  \\
E &  cad(dc^{-1})^2  & d^{-1} a^{-1} c^{-1} (d c a)^2 d\\
 &   &  \\
I &  dc^{-1}d^{-1}a^{-1}d^{-2}a^{-1}c^{-1} & ((d c a)^2 d)^2 c a d c a d
 \\
 &   &  \\
O  & cd^{-2}a^{-1}bc^{-1} &  d^{-1} a^{-1} c^{-1} a b \\
 &   &  \\
U &  (ad^2c^{-1})^3ad^2b^{-1} & d^{-1} a^{-1} c^4 a d b^{-1} d \\
 &   &  \\
T &  (bc^{-1})^2bd^{-2}a^{-1}cd^{-2}a^{-1} & (d^{-1} b)^3 d^{-1} a^{-1} c^{-2} a d \\
 &   &  \\
M &  b(ad^2)^3c^{-1} & d^{-1} b (d^{-1} a^{-1} c^{-1})^3\\
 &   &  \\
L & bd^{-2}a^{-1}c^{-1}bc^{-1}ad^2b^{-1} &    d^{-1} b c a b d^{-1} a^{-1} c a d b^{-1} d  \\
 &   &  \\
K  & c(ad^2)^2c^{-1}bc^{-1}ad^2c^{-1} &  (d^{-1} a^{-1} c^{-1})^2 d^{-1} b d^{-1} a^{-1} c a d \\
 &   &  \\
Y  & c(ad^2)^2c^{-1}dc^{-1}ad^2b^{-1} & (d^{-1} a^{-1} c^{-1})^2 d c a d c^2 a d b^{-1} d \\
 &   &  \\
B  & cd^{-2}a^{-1}(dc^{-1})^2d^{-1}a^{-1}c^{-1} & d^{-1} a^{-1} c^{-1} (a d^2 c a d c)^3 a d c a d \\
 &   &  \\
N & (ad^2c^{-1})^2d(c^{-1}b)^2d^{-1} & d^{-1} a^{-1} c^2 a d (d c a)^2 b d^{-1} b (d^{-1} a^{-1} c^{-1})^2 d^{-1} \\
 &   &  \\
\end{array}
\]
\end{scriptsize}
\end{table}

\vspace{-1.5cm}

\begin{table}[H]
\caption{Correspondence: plaintext alphabet to ciphertext alphabet II}\label{ZuordnungCRyptoFII}
\begin{scriptsize}\[
\begin{array}{c||c|c}
 &   &  \\
& \  U_{f_{\mathrm x_3}} & \ U_{f_{\mathrm x_4}}   \ \\
 &   &  \\
\hline
\hline
 &   &  \\
A &  a^{-1}d^{-1}(ab^3)^2 & b^{-1}a^3d(aca^2)^2\\
 &   &  \\
E & c^{-1}daba^{-1}d^{-1}& ca^2db^{-1}a^3d\\
 &   &  \\
I &  (ba^{-1}d^{-1})^2(a^{-1}d^{-1}c)^2 & (db^{-1}a^3d)^2a^{-2}c^{-1}a^{-2}c^{-1}\\
 &   &  \\
O  & b^{-3}a^{-2}d^{-1} & a^{-2}c^{-1}a^{-1}b^{-1}a^3d\\
 &   &  \\
U &   (ab^3)^4da & (aca^2)^4d^{-1}a^{-3}b\\
 &   &  \\
T &  (a^{-1}d^{-1})^3(b^{-3}a^{-1})^2& (b^{-1}a^3d)^3a^{-2}c^{-1}a^{-3}c^{-1}a^{-1}\\
 &   &  \\
M &   a^{-1}d^{-1}(c^{-1}da)^3 &   b^{-1}a^3d(ca^2)^3 \\
 &   &  \\
L & (a^{-1}d^{-1})^2ca^{-1}d^{-1}ab^3da & b^{-1}a^3da^{-2}c^{-1}b^{-1}a^3daca^2d^{-1}a^{-3}b\\
 &   &  \\
K &  c^{-1}dac^{-1}ab^3 & (ca^2)^2b^{-1}a^3daca^2\\
 &   &  \\
Y  &(c^{-1}da)^2ba^{-1}d^{-1}ab^3da & (ca^2)^2db^{-1}a^3daca^2d^{-1}a^{-3}b\\
 &   &  \\
B  & b^{-3}a^{-1}(ba^{-1}d^{-1})^3a^{-1}d^{-1}c & a^{-2}c^{-1}a^{-1}(db^{-1}a^3d)^3a^{-2}c^{-1}\\
 &   &  \\
N & (ab^3)^2b(a^{-1}d^{-1})^2b^{-1}& aca^3c(a^2db^{-1}a)^2a^2\\
 &   &  \\
\end{array}
\]
\end{scriptsize}
\end{table}

\begin{table}[H]
\caption{Correspondence: plaintext alphabet to ciphertext alphabet III}\label{ZuordnungCRyptoFIII}
\begin{scriptsize}\[
\begin{array}{c||c|c}
 &   &  \\
&  \ U_{f_{\mathrm x_5}}  \ & \ U_{f_{\mathrm x_6}}  \ \\
 &   &  \\
\hline
\hline
 &   &  \\
A  & b^{-1}dc^{-2}b^{-1}a^{-2}b & (c^{-1} b^{-1} a^{-1})^2 c^{-1} b^{-1} \\
 &  &  \\
E  &  cb^{-1}a^{-1}bdc^{-2} & c a^{-1} d c^{-1} b^{-1} \\
 &   &  \\
I  & dc^{-2}dc^{-1}(c^{-1}b^{-1}ab)^2c^{-1} &  (d c^{-1} b^{-1})^2 a c^{-1} a c^{-1} \\
 &   &  \\
O  & b^{-1}adc^{-2} & b c a c^{-1} b^{-1} \\
 &   &  \\
U  & b^{-1}a^{-4}bc^2d^{-1}b & (a^{-1} c^{-1} b^{-1})^3 a^{-1} \\
 &   &  \\
T  & (b^{-1}dc^{-2})^3b^{-1}a^2b & (c^{-1} b^{-1})^2 a b c a\\
 &   &  \\
M  & b^{-1}dc^{-1}(b^{-1}a^{-1}bc)^2b^{-1}a^{-1}b & c^{-1} b^{-1} (c a^{-1})^3 \\
 &   &  \\
L  & b^{-1}dc^{-2}b^{-1}abc^{-1}b^{-1}dc^{-2}b^{-1}a^{-1}bc^2d^{-1}b & c^{-1} b^{-1} a c^{-2} b^{-1} a^{-1} \\
 &   &  \\
K  & cb^{-1}a^{-1}bcb^{-1}a^{-1}dc^{-2}b^{-1}a^{-1}b & c a^{-1} c (a^{-1} c^{-1} b^{-1})^2 \\
 &   &  \\
Y  & (cb^{-1}a^{-1}b)^2dc^{-2}b^{-1}a^{-1}bc^2d^{-1}b & (c a^{-1})^2 d c^{-1} b^{-1} a^{-1} \\
 &   &  \\
B  & b^{-1}ab(dc^{-2})^3b^{-1}abc^{-1} & b c a (d c^{-1} b^{-1})^3 a c^{-1} \\
 &   &  \\
N  & b^{-1}a^{-2}b(dc^{-2}b^{-1})^2 & (a^{-1} c^{-1} b^{-1})^2 d (c^{-1} b^{-1})^2 d^{-1}\\
 &   &  \\
\end{array}
\]
\end{scriptsize}
\end{table}
\vspace{-1.5cm}
\begin{table}[H]
\caption{Correspondence: plaintext alphabet to ciphertext alphabet IV}\label{ZuordnungCRyptoFIV}
\begin{scriptsize}\[
\begin{array}{c||c|c}
 &   &  \\
& \ U_{f_{\mathrm x_7}}  \ & \ U_{f_{\mathrm x_8}}  \\
 &   &  \\
\hline
\hline 
 &   &  \\
A & a^{-3} b^{-1} d c^{-3} a^{-1} (b a^2)^2 a
 & b^{-1} d^{-1} a^{-1}\\
 &   &  \\
E & c^{-1} a^{-1} b a^3 d c^{-3} &  d^{-2} c d^{-1} b^{-1} a d \\
 &   &  \\
I & d c^{-3} d c^{-3} (a^{-3} b^{-1} a c)^2 &  d^{-1} (b^{-1} a)^2 (d c^{-1} d)^2 d \\
 &   &  \\
O & a^{-1} (a^{-2} b^{-1})^2 d c^{-3} &  a d b^{-1} a d\\
 &   &  \\
U & a^{-1} (b a^2)^4 a c^3 d^{-1} b a^3 &  (d^{-1} a^{-1})^5 b\\
 &   &  \\
T & (a^{-3} b^{-1} d c^{-3})^3 a^{-1} (a^{-2} b^{-1})^2 a & (b^{-1} a d)^3 a d a d \\
 &   &  \\
M & a^{-3} b^{-1} d c^{-3} (c^{-1} a^{-1} b a^3)^3 &  b^{-1} a (d^{-1} c d^{-1})^2 d^{-1} c\\
    &  &\\
L &  a^{-3} b^{-1} d c^{-3} a^{-3} b^{-1} a c a^{-3} b^{-1} d c^{-3} a^{-1} b a^3 c^3 d^{-1} b a^3 & b^{-1} a d c^{-1} d^2 b^{-1} d^{-1} a^{-1} b \\
 &   &  \\
K & c^{-1} a^{-1} b a^3 c^{-1} a^{-1} d c^{-3} a^{-1} b a^3 &  (d^{-2} c)^2 b^{-1} \\
 &    &\\
Y & (c^{-1} a^{-1} b a^3)^2 d c^{-3} a^{-1} b a^3 c^3 d^{-1} b a^3 &  (d^{-2} c)^2 d^{-1} b^{-1} d^{-1} a^{-1} b\\
 &    &\\
B & a^{-3} b^{-1} a (d c^{-3})^3 a^{-3} b^{-1} a c &  (a b^{-1})^3 a d c^{-1} d^2\\
 &     &\\
N & a^{-1} (b a^2)^2 a (d c^{-3} a^{-3} b^{-1})^2 & (d^{-1} a^{-1})^2 d^{-1} (b^{-1} a d)^2 d\\
   &  &\\
\end{array}
\]
\end{scriptsize}
\end{table}
\item[3.2.] With these tables he is able to generate the plaintext from Alice, it is 
\begin{align*}
S=& f^{-1}_{\mathrm x_1}\left( dc^{-1}d^{-1}a^{-1}d^{-2}a^{-1}c^{-1} \right) f^{-1}_{\mathrm x_2}\left(    d^{-1}bcabd^{-1}a^{-1}cadb^{-1}d \right) \\ & f^{-1}_{\mathrm x_3}((ba^{-1}d^{-1})^2(a^{-1}d^{-1}c)^2)
 f^{-1}_{\mathrm x_4}\left((ca^2)^2b^{-1}a^3daca^2\right)\\& f^{-1}_{\mathrm x_5}\left( cb^{-1}a^{-1}bdc^{-2} \right)f^{-1}_{\mathrm x_6}\left( bca(dc^{-1}b^{-1})^3ac^{-1} \right)  
 \\& f^{-1}_{\mathrm x_7}\left( a^{-1} (a^{-2} b^{-1})^2 d c^{-3} \right)  f^{-1}_{\mathrm x_8}\left( (ab^{-1})^3adc^{-1}d^2 \right)\\
 =& \texttt{I LIKE BOB}.
\end{align*}
\end{enumerate}
\end{enumerate}
\end{example}


\begin{center}\underline{Security}\end{center}

This private key cryptosystem is secure against chosen plaintext attacks and chosen ciphertext attacks. In a chosen plaintext attack, an attacker, Eve, chooses an arbitrary plaintext of her choice and gets the corresponding ciphertext. In a chosen ciphertext attack Eve sees ciphertexts and gets to some of these ciphertexts  the corresponding plaintexts (see also \cite{BFKR}).\\

An eavesdropper, Eve, intercepts the ciphertext
\begin{align*}
C=c_1  c_2 \cdots  c_z,
\end{align*}
with  $c_i = f_{\mathrm x_i}(u_j)$ for some $1 \leq j \leq N$.
If Alice and Bob choose non characteristic subgroups, then it is likely that $c_j \notin F_U$ for some $1~\leq~j~\leq~z$. Hence the ciphertext units give no hint for the subgroup $F_U$.
Eve knows $L~=~\sum_{k=1}^{z} |c_k|$, the length of $C$, because Alice and Bob are doing no cancellations between $c_i$ and $c_{i+1}$, for $1 \leq i \leq z-1$.

To break the system Eve needs to know the set $U$.
For this it is likely that she assumes that the ball $B(F,L)$ in the Cayleygraph for $F$ contains a basis for $F_U$.
 With this assumption she searches  for primitive elements for $F_U$ in the ball $B(F,L)$, $|y|\leq L$, $y \in F$. In fact she needs to find $N$  primitive elements for $F_U$ in $B(F,L)$ (these would be primitive elements for $F_U$ in a ball $B(F_U,L)$ for some Nielsen reduced basis for $F_U$).
From Proposition \ref{primnumb1} and  Theorem \ref{primnumb2} it is known that the number of primitive elements grows exponentially with the free length of the elements.
Eve chooses sets $M_i~:=~\{m_{i_1},m_{i_2},\ldots,m_{i_K}\}$ with $K\geq N$ and elements $m_{i_j}$ in $B(F,L)$ and with  Nielsen transformations she constructs the corresponding Nielsen reduced sets $M'_i$. If $|M'_i|=N$ then $M'_i$ is a candidate for $U$.\\

 The number $N$ is a constant in the cryptosystem, hence it takes $\mathcal O(\lambda^2)$ time, with $\lambda:=\textnormal{max}\{|m_{j_{\ell}}|  \mid  \ell=1,2,\ldots,K \} \leq L$, to get the set $M'_j$ from $M_j$ with the algorithm \cite{S} (see Remark \ref{AlgNred}).

The main security certification depends on the fact, that for a single subset of $K$ elements Eve finds a Nielsen reduced set in polynomial running time (more precisely in quadratic time) but she has to test all possible subsets of $K$ elements for which she needs exponential running time.\\

The security certification can be improved by the next two improvements.\\

First, Alice and Bob choose in addition an explicit presentation of the ciphertext units $c_i$ as matrices in SL$(2,\mathbb Q)$.
So, they agree on a faithful representation
\begin{align*}
 \varphi: F &\rightarrow SL(2, \mathbb Q)\\
 x_i &\mapsto M_i,
\end{align*}
 of $F$ into $SL(2,\mathbb Q)$ (see Theorem \ref{TL}).
The group $G=\varphi(F)$ is isomorphic to $F$ under the mapping   $x_i \mapsto M_i$, for $i=1,\ldots,q$. The ciphertext is now 
\begin{align*}
C'&= \varphi(c_1) \varphi(c_2) \cdots \varphi(c_z)\\
 &=W_1W_2\cdots W_z,
\end{align*}
a sequence of matrices $W_j \in \text{SL}(2,\mathbb Q)$. The encryption is realizable with a table (as Table \ref{ZuordnungCRyptoFabstract}) if the representation $\varphi$ is applied to the elements in the table. Therefore Bob gets a table with matrices and hence an assignment from the matrices to the plaintext alphabet depending on the automorphisms  $f_{\mathrm x_i}$.

Here the additional security certification is, that there is no algorithm known to solve the  membership problem (see for instance \cite{MKS}) for subgroups of $\text{SL}(2,\mathbb Q)$ which are not subgroups in  SL$(2,\mathbb Z)$. 
 B. Eick, M. Kirschner and C. Leedham-Green presented in the paper \cite{EKLG} a practical algorithm to solve the constructive membership problem for discrete free subgroups of rank $2$ of SL$(2, \mathbb R)$. For example, the subgroup SL$(2,\mathbb Z)$ of  SL$(2, \mathbb R)$ is discrete.
But they also mention, that it is an open problem to solve the membership problem for  arbitrary subgroups of SL$(2,\mathbb R)$ with rank $m \geq 2$. Alice and Bob work with  subgroups of rank $N \geq 2$. Hence there is in general no algorithm known for Eve to solve the membership problem, in particular there is always no such algorithm known for  $N\geq 3$.

\begin{example}
In this example\footnote{We realized this example with the computer programs Classic Worksheet Maple~16 and GAP. In GAP we used the package ``FGA'' (Free Group Algorithms, a GAP4 Package by Christian Sievers, TU Braunschweig).}  Alice and Bob agree additionally to Example \ref{example1} on a faithful representation. With Theorem \ref{TL} they generate the matrices 
\begin{align*}
X_1:=\left( \begin{smallmatrix}  \frac{-7}{2} & \frac{45}{4}\\
                                                1 & \frac{-7}{2} \end{smallmatrix} \right), \
X_2:=\left( \begin{smallmatrix}  \frac{-15}{2} & \frac{221}{4}\\
                                                1 & \frac{-15}{2} \end{smallmatrix} \right) \text{ and }
X_3:=\left( \begin{smallmatrix}  \frac{-23}{2} & \frac{525}{4}\\
                                                1 & \frac{-23}{2} \end{smallmatrix} \right).
\end{align*}
These matrices form a basis for a free group $G$ of rank $3$. Alice and Bob generate a subgroup $G_1$ of $G$ with rank $4$. The free generating set for $G_1$ is $\{X_1X_2, X_3X_1^2,X_2X_3X_2,X_1^{-1}X_2\}$.
They choose the faithful representation 
\begin{align*}
\varphi: F &\to SL(2,\mathbb Q)\\
         a & \mapsto X_1X_2= \left( \begin{smallmatrix} \frac {75}{2}  & \frac {-1111}{4}  \\ -11 & \frac {163}{2} \end{smallmatrix} \right),
		&b & \mapsto X_3X_1^2= \left( \begin{smallmatrix} -1189 & 3990 \\
104 & -349 \end{smallmatrix} \right),\\
        c & \mapsto X_2X_3X_2= \left( \begin{smallmatrix} -2681 & 19966 \\
360 & -2681 \end{smallmatrix} \right),
       &d & \mapsto X_1^{-1}X_2=\left( \begin{smallmatrix} 15 & -109 \\
4 & -29 \end{smallmatrix} \right).
\end{align*}
The ciphertext is now 
\begin{align*}
C'=& \varphi( dc^{-1}d^{-1}a^{-1}d^{-2}a^{-1}c^{-1})   \ \varphi( d^{-1}bcabd^{-1}a^{-1}cadb^{-1}d) \\  &\varphi( (ba^{-1}d^{-1})^2(a^{-1}d^{-1}c)^2)  \  \varphi( (ca^2)^2b^{-1}a^3daca^2) \  \varphi( cb^{-1}a^{-1}bdc^{-2})   \\ 
& \varphi( bca(dc^{-1}b^{-1})^3ac^{-1} )   \
 \varphi( a^{-1} (a^{-2} b^{-1})^2 d c^{-3}) \   \varphi( (ab^{-1})^3adc^{-1}d^2)\\
=& \left( \begin{smallmatrix} \frac {-429743093559909}{2}  & \frac {-6400784021410159}{4}  \\ 
-62588240305379 &  \frac {-932216979117085}{2}   \end{smallmatrix} \right)\\
&  \left( \begin{smallmatrix}   \frac {-3240070331754423030683243991}{2}  & 
 \frac {47007695458416827592369656315}{4}  \\  -223326322203710575272321977 & \frac {
3240070327830150751386194361}{2} \end{smallmatrix} \right)\\
&\left( \begin{smallmatrix} \frac {-6899014060703475554169965}{2}  &  \frac {102756972145191520348785607}{4}  \\ 
301722468685102729969483 &  \frac {-4493988131847945704997109}{2}    \end{smallmatrix} \right)\\
&
\left( \begin{smallmatrix}   \frac {-397074726172421275253684843812134445}{2} & \frac{5883318761059670223751985896578473377}{4} \\  26659253089426526822952736194350493  &
\frac {-395000924306510751052288425218790757}{2}  \end{smallmatrix} \right)\\
&\left( \begin{smallmatrix}  \frac {46475888407425825}{2}  & \frac {692232489736400389}{4}  \\ 
-3120351373297111 & \frac {-46475896943687759}{2} \end{smallmatrix} \right)\\
& \left( \begin{smallmatrix} \frac {-37154085868492177463035768197599}{2}  & 
 \frac {-553374013794643763898030444104547}{4}  \\ 
1624906569753714749910956723073 & \frac {24201404758781402065719318991873}{2} \end{smallmatrix} \right)\\
&\left( \begin{smallmatrix}  \frac {-3418963163764785449276501363}{2} & \frac {-50923553357916815212095363641}{4}  \\
-230751369629481141540301125 & \frac {-3436913216344813651054341083}{2} \end{smallmatrix} \right)\\
& \left( \begin{smallmatrix} \frac {2739747352948144349387}{2}  & \frac {-39628644296581967709615}{4} \\ 
-402070084312200114547 &  \frac {5815679440792026855107}{2}\end{smallmatrix}\right).
\end{align*}
Instead of a sequence of words in $F$ Alice sends to Bob a sequence of eight matrices in SL$(2,\mathbb Q)$.
\end{example}

For the second improvement Alice and Bob use instead of a presentation of the ciphertext in SL$(2,\mathbb Q)$ a presentation of the ciphertext in a free group in GL$(2, k)$ with $k~:=~\mathbb Z[y_1,y_2,\ldots,y_w]$, the ring of polynomials in variables $y_1,y_2,\ldots,y_w$. With the help of a homomorphism $\epsilon ^*: \text{GL}(2,k) \to \text{GL}(2,\mathbb Z)$ and the knowledge of an algorithm to write each element in the modular group PSL$(2,\mathbb Z)$, the group of $2 \times 2$ projective integral matrices of determinant $1$, in terms of $s$ and $t$ they can reconstruct the message. Here,  
\begin{align*}
s=\left( \begin{smallmatrix}
  \phantom{-}0 & 1 \\
 -1 & 0
 \end{smallmatrix} \right) \quad \text{and} \quad 
t=
\left( \begin{smallmatrix}
 1 & 1 \\ 0 & 1
 \end{smallmatrix} \right)
\end{align*}
and  $\text{PSL}(2,\mathbb Z) = \langle s, t \mid s^2=(st)^3 = 1  \rangle$.

Every finitely generated free group is faithfully represented by a subgroup of the modular group PSL$(2,\mathbb Z)$. 
Especially, the two matrices
\begin{align*}
\left( \begin{smallmatrix}
  \phantom{-}0 & 1 \\
 -1 & 2
 \end{smallmatrix} \right) \quad \text{and} \quad 
\left( \begin{smallmatrix}
 \phantom{-}2 & 1 \\ -1 & 0
 \end{smallmatrix} \right)
\end{align*}
generate a free group of rank two, and this free group
certainly contains finitely generated free groups.

This improvement is very similar to the version in \cite{BFKR}. 
Here, the security certification depends in addition on the unsolvability of Hilbert's Tenth Problem.  Y. Matiyasevich proved in \cite{Ma} that there is no general algorithm which determines whether or not an integral polynomial in any number of variables has a zero.

\section{Public Key Cryptosystem Based on Automorphisms of Free Groups $F$}\label{ElGamal}

Now we describe a public key cryptosystem for Alice and Bob which is inspired by the ElGamal cryptosystem (see \cite{E} or \cite[Section 1.3]{MSU}), based on discrete logarithms, that is:
\begin{enumerate}
\item Alice and Bob agree on a finite cyclic group $G$ and a generating element $g \in G$.
\item Alice picks a random natural number $a$ and publishes the element $c:=g^a$.
\item Bob, who wants to transmit a message $m \in G$ to Alice, picks a random natural number $b$ and sends the two elements $m\cdot c^b$ and $g^b$, to Alice. Note that $c^b=g^{ab}$.
\item Alice recovers $m=\left( m \cdot c^b\right) \cdot \left( \left( g^b \right)^a\right)^{-1}$. 
\end{enumerate}

Let $X=\{x_1,x_2,\ldots,x_N\}$, $N\geq 2$, be the free generating set of the free group $F=\langle X \mid  \phantom{R} \rangle$. It is $X^{\pm 1}=X \cup X^{-1}$.
The message is an element $m \in S^*$, the set of all freely reduced words with letters in $X^{\pm 1}$.
Public are the free group $F$, its free generating set $X$ and an element $a \in S^*$.
 The automorphism $f$ should be chosen randomly, for example as it is described in Section \ref{chooseNT}.\\

The public key  cryptosystem is now as follows:\\

\textbf{Public parameters:} The group $F= \langle X \mid \phantom{R} \rangle$, a freely reduced word  $ a \ne 1$ in the free group $F$ and an automorphism $f: F \to F$ of infinite order.\\

\textbf{Encryption and Decryption Procedure:}
\begin{enumerate}
\item Alice chooses privately a natural number $n$ and publishes the element $f^n(a)=:c \in S^*$.

\item Bob picks privately a random $t \in \mathbb N$ and his message $m\in S^*$.
He calculates the freely reduced elements
\begin{align*}
m \cdot f^t(c)=:c_1 \in S^* \ \text{ and } \ f^t(a)=:c_2 \in S^*.
\end{align*} He sends  the ciphertext $(c_1,c_2) \in S^* \times S^*$ to Alice. 

\item Alice calculates
\begin{align*}
c_1 \cdot f^n(c_2)^{-1}&= m\cdot f^t(c) \cdot f^n(c_2)^{-1}\\
                          &= m\cdot f^t(f^n(a)) \cdot (f^n(f^t(a))^{-1}\\
                          &= m\cdot f^{t+n}(a) \cdot (f^{n+t}(a))^{-1}\\
                          &= m,
\end{align*}
and gets the message $m$.
\end{enumerate}

\begin{remark}
A possible attacker, Eve,  can see the elements  $c, c_1, c_2 \in S^*$. She does not know the free length of $m$ and the cancellations between $m$ and $f^t(c)$ in  $c_1$.
 It could be possible that $m$ is completely canceled by the first letters of $f^t(c)$.
 Hence she cannot determine  $m$ from the  given  $c_1$.   Eve just sees words, $f^t(a)$ and $f^n(a)$, in the free generating set $X$ from which it is unlikely to realize the exponents $n$ and $t$, that is, the private keys from Alice and Bob, respectively. The security certification is  based on the Diffie-Hellman-Problem.
\end{remark}

\begin{remark} We give some ideas to  enhance the security, they can also be combined:
\begin{enumerate}
\item The element $a \in S^*$ could be taken as a common private secret between Alice and Bob. They could use for example the Anshel-Anshel-Goldfeld key exchange protocol (see \cite{MSU}) to agree on the element $a$. 
\item Alice and Bob agree on a faithful representation from $F$ into the special linear group of all $2 \times 2$ matrices with entries in $\mathbb Q$, that is, $g: F \to \text{SL}(2, \mathbb Q)$. Now $m \in S^*$ and Bob sends $g(m)  \cdot g(f^t(c))=:c_1 \in \text{SL}(2, \mathbb Q)$ instead of $m \cdot f^t(c)=:c_1 \in S^*$; $c$ and $c_2$ remain the same. Therefore, Alice calculates $c_1 \cdot (g(f^n(c_2)))^{-1}=g(m)$ and hence the message $m=g^{-1}(g(m))\in S^*$. 
This variation in addition extends the security certification to the membership problem in the matrix group $\text{SL}(2, \mathbb Q)$ (see \cite{EKLG}).
\end{enumerate}    
\end{remark}

\begin{example}~\\
This example\footnote{We used the computer program GAP and the package ``FGA'' (Free Group Algorithms, a GAP4 Package by Christian Sievers, TU Braunschweig).} is a very small one and it is just given for illustration purposes. Bob wants to send a message to Alice. \\
 The \textbf{public parameters} are the free group $F=\langle x_1,x_2,x_3 \mid \phantom{R} \rangle$ of rank $3$, the freely reduced word $a \in F$, with $a:=x_1^2 x_2 x_3^{-2} x_2$ and the  automorphism $f:F \to F$, which is given, for this example,  by regular Nielsen transformations:\\
$f=[(N2)_{1.2}]^2 \ (N2)_{3.2} \ (N1)_3 \ (N2)_{2.3} $, that is,
\begin{align*}  
x_1 \mapsto x_1x_2^2, \
x_2 \mapsto x_3^{-1}, \
x_3 \mapsto x_2^{-1}  x_3^{-1}.
\end{align*}
\begin{enumerate}
\item Alice's private key is $n=7$. Thus, she  gets the automorphism \begin{align*}
f^7:F &\to F\\   
x_1 &\mapsto x_1x_2^2x_3^{-1}x_2(x_2x_3)^2(x_3x_2x_3^2x_2)^2x_3x_2\\
x_2 &\mapsto x_2^{-1}((x_3^{-1}x_2^{-1}x_3^{-1})^2x_2^{-1}x_3^{-1})^2x_3^{-1}x_2^{-1}x_3^{-2}\\
x_3 &\mapsto (((x_2^{-1}x_3^{-1})^2x_3^{-1})^2x_2^{-1}x_3^{-2})^2x_2^{-1}(x_3^{-1}x_2^{-1}x_3^{-1})^2x_3^{-1}.
\end{align*}
Her public key is 
\begin{align*}
c:=f^{7}(a)=&(x_1x_2^2x_3^{-1}x_2(x_2x_3)^2(x_3x_2x_3^2x_2)^2x_3x_2)^2(x_3^2x_2)^2\\&((x_3x_2x_3)^2x_2x_3)^2x_3x_2x_3^2x_2x_3^{-1}.
\end{align*}
\item Bob privately picks  the ephemeral key $t=5$ and gets the automorphism 
\begin{align*}
f^5:F &\to F\\   
x_1 &\mapsto x_1x_2^2x_3^{-1}x_2^2x_3(x_3x_2)^2\\
x_2 &\mapsto x_2^{-1}(x_3^{-1}x_2^{-1}x_3^{-1})^2x_3^{-1} \\
x_3 &\mapsto ((x_2^{-1}x_3^{-1})^2x_3^{-1})^2x_2^{-1}x_3^{-2}.
\end{align*}
His message for Alice is $m=x_3^{-2}x_2^2x_3x_1^2x_2^{-1}x_1^{-1}$.
He calculates 
\begin{align*}
c_1=& m \cdot f^5(c)\\
=& x_3^{-2} x_2^2 x_3 x_1^2 (x_2 x_3^{-1})^2 ((x_3^{-1} x_2^{-1} x_3^{-2} x_2^{-1})^2 x_3^{-2} x_2^{-1})^2 (x_3^{-1} x_2^{-1} x_3^{-1})^2 x_3^{-1} x_2^{-1}\\ 
&((((x_3^{-1} x_2^{-1} x_3^{-1})^2 x_2^{-1} x_3^{-1})^2 x_3^{-1} x_2^{-1} x_3^{-1} x_2^{-1} x_3^{-1})^2 (x_3^{-1} x_2^{-1} x_3^{-2} x_2^{-1})^2 x_3^{-1}\\
& x_2^{-1} x_3^{-1})^2 ((x_3^{-1} x_2^{-1} x_3^{-2} x_2^{-1})^2 x_3^{-2} x_2^{-1})^2 (x_3^{-1} x_2^{-1} x_3^{-1})^2 x_3^{-1} x_1 x_2^2 x_3^{-1} x_2 (x_3^{-1}\\
& (((x_3^{-1} x_2^{-1} x_3^{-2} x_2^{-1})^2 x_3^{-2} x_2^{-1})^2 (x_3^{-1} 
x_2^{-1} x_3^{-1})^2 x_3^{-1} x_2^{-1})^3 (x_3^{-1} x_2^{-1} x_3^{-1})^2 \\
&x_2^{-1} x_3^{-1} ((x_3^{-1} x_2^{-1} x_3^{-2} x_2^{-1})^2 x_3^{-2} x_2^{-1})^2 (x_3^{-1} x_2^{-1} x_3^{-1})^2 x_2^{-1})^3 x_3^{-1}\\
& ((x_3^{-1} x_2^{-1} x_3^{-2} x_2^{-1})^2 x_3^{-2} x_2^{-1})^2 (x_3^{-1} x_2^{-1} x_3^{-1})^2 x_2^{-1} x_3^{-1} x_2
\end{align*}
and 
\begin{align*}
c_2:=f^5(a)=(x_1x_2^2x_3^{-1}x_2^2x_3(x_3x_2)^2)^2x_3^2x_2(x_3x_2x_3)^2x_3x_2x_3^{-1}.
\end{align*}
The ciphertext for Alice is the tuple $(c_1,c_2)$.
\item 
Alice first computes
\begin{align*}
(f^7(c_2))^{-1}=&x_2^{-1} (((((x_3 x_2)^2 x_3)^2 x_3 x_2 x_3)^2 x_3 x_2 (x_3 x_2 x_3)^2)^2 x_3 x_2 \\& ((x_3 x_2 x_3)^2 x_2 x_3)^2 x_3 x_2 x_3)^2  (x_3 x_2 (((x_3 x_2 x_3)^2 x_2 x_3)^2 x_3 x_2 x_3 x_2 x_3)^2 \\& (x_3
 x_2 x_3^2 x_2)^2 x_3)^2 x_2 (((((x_3^2 x_2)^2 x_3 x_2)^2 x_3^2 x_2 x_3 x_2)^2 x_3 \\& (x_3 x_2 x_3^2 x_2)^2 x_3 x_2)^2 x_3 (x_3 x_2 x_3^2 x_2)^2 x_3^2 x_2\\
& (((x_3 x_2 x_3)^2 x_2 x_3)^2 x_3
 x_2 x_3 x_2 x_3)^2 (x_3 x_2 x_3^2 x_2)^2 x_3 \\& (x_3 x_2^{-1})^2 x_2^{-1} x_1^{-1})^2
\end{align*}
and gets $m$ by
$$
m=c_1 \cdot (f^7(c_2))^{-1}=x_3^{-2}x_2^2x_3x_1^2x_2^{-1}x_1^{-1}.
$$
\end{enumerate}
\end{example}

\section{Conclusion}
In comparison to the standard cryptosystems which are mostly based on number theory we explained two cryptosystems which use combinatorial group theory.
The first cryptosystem in Section \ref{privkeycryp1} is a one-time pad, which choice of the random sequence for encryption is not number-theoretic. At the moment it is costlier than the standard systems but it is another option for a one-time pad which is based on combinatorial group theory and not on number theory.
The second cryptosystem in Section \ref{ElGamal} is similar to the ElGamal cryptosystem  (see~\cite{E}), which is easier to handle. The ElGamal cryptosystem is based on the discrete logarithm problem over a finite field. If this problem should  eventually be solved we introduced here an alternative system, which is not based on number theory.

\newpage

\bibliographystyle{abbrv}

\end{document}